\documentclass[preprint, a4paper, 11pt, oneside, onecolumn]{elsarticle}
\usepackage{amssymb, amsthm}
\usepackage{subcaption,caption}
\usepackage{graphicx}
\usepackage{float}    
\usepackage{amsthm}
\usepackage{amsmath}
\usepackage{mathrsfs}
\usepackage[misc]{ifsym}
\usepackage{diagbox}
\usepackage{graphicx} 
\usepackage{amsmath,latexsym,amssymb,amsfonts,amsbsy}
\usepackage{algorithm}
\usepackage{algpseudocode}
\usepackage{color,xcolor}
\usepackage{booktabs}
\usepackage{threeparttable}
\usepackage{multicol}
\usepackage{multirow}
\usepackage{epstopdf}
\numberwithin{equation}{section}
\newtheorem{theorem}{Theorem}[section]

\newtheorem{proposition}[theorem]{Proposition}

\begin{document} 
\begin{sloppypar}
    \begin{frontmatter}
       \title{Adaptive Parameter Selection in Nudging Based Data Assimilation}
\author[1]{Aytekin \c{C}{\i}b{\i}k}
\ead{abayram@gazi.edu.tr}

\author[2]{Rui Fang}
\ead{ruf10@pitt.edu}
\author[2]{William Layton\corref{cor1}}
\ead{wjl@pitt.edu}
\cortext[cor1]{Corresponding author: William Layton}

\author[2]{Farjana Siddiqua}
\ead{fas41@pitt.edu}
\affiliation[1]{organization={Department of Mathematics, Gazi University},
        city={Ankara},
postcode={06550},
country={Turkiye}}

\affiliation[2]{
  organization={Department of Mathematics, University of Pittsburgh},
  city={Pittsburgh},
  postcode={15260},
  country={USA}
}

\begin{abstract}
Data assimilation combines (imperfect) knowledge of a flow's physical laws with (noisy, time-lagged, and otherwise imperfect) observations to
produce a more accurate prediction of flow statistics. Assimilation by nudging (from 1964), while non-optimal, is easy to implement and its analysis is clear and well-established. Nudging's uniform in time accuracy has even been established under conditions on the nudging parameter $\chi$ and the density of observational locations, $H$,  Larios, Rebholz, and Zerfas \cite{larios2019global}. One remaining issue is that nudging requires the user to select a key parameter. The conditions required for this parameter, derived through $\acute{a}$ priori (worst case) analysis are severe (Section \ref{aprior-analysis}  herein) and far beyond those found to be effective in computational experience. One resolution, developed herein, is self-adaptive parameter selection. This report develops, analyzes, tests, and compares two methods of self-adaptation of nudging parameters. One combines analysis and response
to local flow behavior. The other is based only on response to flow behavior. The comparison finds both are easily implemented and yield effective values of the nudging parameter much smaller than those of $\acute{a}$ priori analysis.
\end{abstract}

\begin{keyword}nudging, adaptive algorithm, FEM, data assimilation\\
        \vspace{0.2cm}
        \end{keyword}   
\end{frontmatter}
\section{Introduction}

Predictions of the future state of a flow, here internal 2d or 3d flow of an
incompressible viscous fluid in a domain $\Omega$,
\begin{align}
u_{t}+u\cdot\nabla u-\nu\triangle u+\nabla p  &  =f(x)\text{, and }\nabla\cdot
u=0,\text{ in }\Omega, \ 0<t\leq T,\\
u  &  =0\text{ on }\partial\Omega\text{ and }u(x,0)=u_{0}(x),
\end{align}
are improved, Kalnay \cite{kalnay2002atmospheric}, by incorporating / assimilating observations
of the flow
\[
u_{obs}(x,t)\text{ for }0<t\leq T.
\]
The goal of data assimilation is to combine incomplete, sparse, noisy, and
possibly time-delayed observations, $u_{obs}$, with incomplete and approximate
knowledge of a flow's dynamic laws to produce a more accurate prediction of
flow statistics. Nudging-based assimilation, from 1964 Luenberger \cite{luenberger1964observing},
is amenable to both analysis, e.g. Biswas and Price \cite{biswas2021continuous}, Azouani,
Olson and Titi \cite{azouani2014continuous}, Cao, Giorgini, Jolly and Pakzad \cite{cao2022continuous},
Larios, Rebholz, and Zerfas \cite{larios2019global} (among hundreds of papers), and straightforward
implementation. However, nudging is non-optimal since no criterion is
minimized, Lakshmivarahan and Lewis \cite{lakshmivarahan2013nudging}. Uniform in-time convergence
has even been proven \cite{larios2019global} under parameter conditions reviewed in Section \ref{continumm-nudging}. A difficulty
remaining is that users must select effective nudging parameters and
the analytical theory places parameter restrictions, see Section \ref{aprior-analysis}, far beyond ranges found effective in computational experience, summarized in Kalnay \cite{kalnay2002atmospheric}. Since \'{a} priori analysis develops conditions from a series of worst-case estimates, these
restrictions may be pessimistic, and self-adaptive parameter selection,
developed herein, is natural.

To develop adaptive-$\chi$ methods, we adopt the simplest interesting setting. Let the true velocity be denoted $u(x,t)$ and its observed values $u_{obs}=I_{H}u$ be an $L^{2}$ projection of the true values on a finite-dimensional subspace (associated
with a length-scale denoted $H$). The continuum, nudged approximation $v(x,t)$
satisfies
\begin{gather}
\text{select parameter }\chi,\text{ set\ }v(x,0)=v_{0}(x)\text{ and
solve}\nonumber\\
v_{t}+v\cdot\nabla v-\nu\triangle v-\chi I_{H}(u-v)+\nabla q=f(x)\text{,
}\nabla\cdot v=0, \text{ in } \Omega.
\end{gather}
The initial and boundary conditions are
\begin{equation*}
\begin{gathered}
u=0 \text{ on } \partial \Omega, \ v=0 \text{ on } \partial \Omega, \
u(x,0) =u_0(x), \text{ and } v(x,0) = v_0(x).
\end{gathered}
\end{equation*}
To support the necessity for adapting $\chi$, in Section \ref{continumm-nudging} we summarized error
analysis (inspired by the analysis in \cite{larios2019global}) showing that for $H$
\textit{small enough} and $\chi$\ \textit{large enough} nudging is accurate uniformly in time. Specifically, the conditions from this analysis are
\begin{align}
\text{ }H-\text{condition}  &  \text{: }\nu-2C_{1}^{2}H^{2}\chi\geq
0,\label{eqResult1}\\
\text{2d }\chi-\text{condition}  &  \text{: }\chi-\frac{1}{2\nu}\frac{1}%
{T}\int_{0}^{T}||\nabla u||^{2}dt\geq\chi_{0}>0,\label{eq2dChi}\\
\text{3d }\chi-\text{condition}  &  \text{: }2\left[  \chi-\frac{2048}%
{19683}\nu^{-3}\frac{1}{T}\int_{0}^{T}||\nabla u||^{4}dt\right]  \geq\chi
_{0}>0. \label{eq3dChi}%
\end{align}
For high Reynolds number flows in 2d and 3d, in Section \ref{aprior-analysis}  we investigate
what turbulence phenomenology suggests \textit{large enough} and \textit{small
enough} in (\ref{eqResult1})-(\ref{eq3dChi}) mean. There results
\begin{align}
\text{in 2d} & \text{: }\chi \gtrsim \mathcal{O}(\mathcal{R}e^{3/2})\text{ \ and
\ }H \lesssim \mathcal{O}(\mathcal{R}e^{-5/4}), \label{eqReBasedInterpretation}\\
\text{in 3d}  &  \text{: }\chi \gtrsim \mathcal{O}(\mathcal{R}e^{5})\text{ \ and
\ \ }H \lesssim \mathcal{O}(\mathcal{R}e^{-3}).\label{3d-condition}
\end{align}

While sufficient (in theory), these are impractical (in practice).
Extensive practical experience, Kalnay \cite{kalnay2002atmospheric}, shows that the nudging parameter is best chosen of intermediate size, balancing data errors with
model errors and numerical errors. Thus the parameter requirements
(\ref{eqResult1})-(\ref{eq3dChi}) are far beyond parameter values in problems
where assimilation is used. The gap is likely because of conditions
(\ref{eqResult1})-(\ref{eq3dChi}) emerge from analysis through a series of
worst-case estimates which may not occur simultaneously and (\ref{eqReBasedInterpretation}), (\ref{3d-condition}) emerge from the assumption that the flow is fully developed turbulent. To address this
issue, in Section \ref{adapt-section} we present 2 low complexity algorithms for self-adaptive
parameter selection. These act to find a smaller effective $\chi$\ value in
response to the local flow behavior and are only slightly more expensive in
operations and memory than \'{a} priori selection. 

Adaptivity exploits $\chi$-sensitivity to improve approximations. It is based on $3$ principles: localization (of the global criteria), estimation (of the deviation from a desired state), and decision (of how to choose the next parameter value). 
The condition that $\|v(\cdot, t)- u(\cdot, t)\|$ is decreasing is localized to decreasing every time step. The estimator is $\|I_H(u-v)\|$ and the decision tree is to increase $\chi$ for large estimates and decrease $\chi$ for very small estimators. We compare two decision strategies based on derived analytical criteria for decrease and simply asking if a reduction is observed.

Section \ref{sec: tests} presents tests of the adaptive nudging methods for problems with model errors with both \'{a} priori and self-adaptive parameter selections. We
find both adaptive-$\chi$ methods produce effective $\chi$ values smaller than (\ref{eqReBasedInterpretation})-(\ref{3d-condition}) and accurate velocities. We also observe that an H-condition seems to be necessary for accurate, long-time approximations. Addressing the H-condition is an open problem. In many applications, $H$ depends on the density of sensors. 

\subsection{Related work}

Nudging for data assimilation stems from the 1964 work of Luenberger \cite{luenberger1964observing}
and has extensive tests in geophysical flow simulations, e.g., \cite{hoke1976initialization, kalnay2002atmospheric, lakshmivarahan2013nudging, navon1998practical, stauffer1993optimal, zou1992optimal}. Its form also arises
under the terms various types of damping and time relaxation, e.g.,
\cite{breckling2017review}, with large numbers of papers for each realization. The
analytical work on nudging is similarly enormous; we mention three that we
(and most others) build upon. Azouani, Olson, and Titi \cite{azouani2014continuous} gave an
early and possibly even the first complete convergence analysis of a nudging
algorithm. Their analysis produced a similar $H$-condition as herein in
Theorem 2.1 p. 254. Biswas and Price \cite{biswas2021continuous} gave a remarkable analysis in
3d without assuming\ \'{a} priori strong solutions. Their analysis recovered
similar conditions as herein on $H$ (their equation (3.27)) and $\chi$\ (their
(3.8)). Biswas and Price \cite{biswas2021continuous} also provided an early algorithm for
adaptation of $\chi$ which retained a Reynolds number dependent lower bound on
$\chi$\ in their equation (4.4). In 2019 Larios, Rebholz, and Zerfas
\cite{larios2019global} proved uniform in time error bounds in a way that narrowed the
gap between analysis and computational tests.

None of these papers sought an optimal $\chi$. Nudging emerged from
Luenberger's work in control theory then developed in independent directions.
Along its path, several reconnected nudging to control and optimization,
determining methods for optimization of the nudging parameter. For work on this approach (not considered herein) involving solving the adjoint problem,
see Hoke and Anthes \cite{hoke1976initialization}, Stauffer and Bao \cite{stauffer1993optimal}, Zou, Navon,
LeDimet \cite{zou1992optimal} and Navon \cite{navon1998practical}.
Carlson, Hudson, and Larios \cite{carlson2020parameter} and Martinez \cite{martinez2022convergence} have studied recovering the unknown viscosity parameter for the 2D Navier-Stokes equations (NSE). Parameter estimates for the Lorenz equation are studied in Carlson, Hudson, Larios, Martinez, Ng, and Whitehead \cite{carlson2021dynamically}.

\subsection{Preliminaries}
The notation for the $L^2(\Omega)$ norm and inner product is $\|\cdot\|$ and $(\cdot,\cdot)$, respectively. Using $\|\cdot\|_{L^p}$, we indicate the $L^p(\Omega)$ norm. 
The numerical tests use a standard finite element discretization, see \cite{layton2008introduction} for details. We also use the inequalities derived by Ladyzhenskaya \cite{ladyzhenskaya1969mathematical}.

\begin{theorem}(The Ladyzhenskaya Inequalities, see Ladyzhenskaya \cite{ladyzhenskaya1969mathematical}) For any vector function $u:\mathbb{R}^d\rightarrow \mathbb{R}^d$ with compact support and with the indicated $L^p$ norms finite,
\begin{equation*}
\begin{gathered}
\|u\|_{L^4(\mathbb{R}^2)}\leq 2^{1/4}\|u\|_{L^2(\mathbb{R}^2)}^{1/2}\|\nabla u\|_{L^2(\mathbb{R}^2)}^{1/2},\ (d=2),\\
\|u\|_{L^4(\mathbb{R}^3)}\leq \frac{4}{3\sqrt{3}}\|u\|^{1/4}\|\nabla u\|^{3/4},\ (d=3),\\
\|u\|_{L^6(\mathbb{R}^3)}\leq \frac{2}{\sqrt{3}}\|\nabla u\|,\ (d=3).
\end{gathered}
\end{equation*}
\end{theorem}

\section{Continuum Nudging}\label{continumm-nudging}
To motivate the necessity of adaptivity, for completeness, we review the (now
standard since \cite{biswas2021continuous}, \cite{azouani2014continuous}, \cite{larios2019global}) proof that the error
$\rightarrow0$\ uniformly in $T$ as $\chi\rightarrow\infty$ and uniformly in
$\chi$\ as $T\rightarrow\infty.$ The (direct) proof shows how conditions
$H$\textit{\ small enough and }$\chi$\textit{\ large enough} emerge
organically for the structure of the NSE and it motivates one
adaptive algorithm of Section \ref{adapt-section}. This conditions $H$\textit{\ small
enough and }$\chi$\textit{\ large enough}\ are interpreted under 2d and 3d
turbulence phenomenology in Section \ref{aprior-analysis} below, leading to
(\ref{eqReBasedInterpretation}). We assume here that $I_{H}$\ is an $L^{2}$
projection and satisfies%
\begin{equation}
||(I-I_{H})w||\leq C_{1}H||\nabla w||\text{ for all }w\in\left(  H_{0}%
^{1}(\Omega)\right)  ^{d},d=2\text{ or }3. \label{eqInterpError}
\end{equation}
For the 3d case, we assume\footnote{One can assume a stronger condition, such
as \ $||\nabla u||\in L^{\infty}(0, T)$, and obtain an estimate that
superficially appears better behaved concerning the Reynolds number or a weaker condition, such as Prodi-Serrin, and obtain an estimate that has a worse explicit $Re$ dependence. This is because some dependence is subsumed within the assumption.} that $||\nabla
u||^{4}\in L^{1}(0, T).$ This assumption is not necessary (since the work of
Biswas and Price \cite{biswas2021continuous}) but shortens the analysis.

\begin{proposition}
Let $e=u-v.$ Suppose $I_{H}$ is an $L^{2}$ projection with the approximation property (\ref{eqInterpError}) above. In 2d let
the parameter conditions (\ref{eqResult1}), (\ref{eq2dChi}) hold and in 3d
(\ref{eqResult1}), (\ref{eq3dChi}). Then, the error $\rightarrow
0$\ uniformly in $T$ as $\chi\rightarrow\infty$ and uniformly in $\chi$\ as
$T\rightarrow\infty$. In particular,
\[
||e(T)||\leq\exp\left[  -\chi_{0}T\right]  ||e(0)||.
\]
\end{proposition}
\begin{proof}

\textbf{The 2d case:} By subtraction, $e$ satisfies
\begin{gather*}
\frac{1}{2}\frac{d}{dt}||e||^{2}+(u\cdot\nabla u-v\cdot\nabla v,e)+\nu||\nabla
e||^{2}+\chi||I_{H}e||^{2}=0, \\
\text{where }|(u\cdot\nabla u-v\cdot\nabla v,e)|=|(e\cdot\nabla u,e)|\leq
||\nabla u||||e||_{L4}^{2}.
\end{gather*}
Note that $||e||^{2}=||I_{H}e||^{2}$ $+||(I-I_{H})e||^{2}$ $\leq||I_{H}%
e||^{2}+$ $\left(  C_{1}H||\nabla e||\right)  ^{2}.$ By the 2d Ladyzhenskaya
and arithmetic-geometric mean inequalities%
\[
|(e\cdot\nabla u,e)|\leq\sqrt{2}||\nabla u||||e||||\nabla e||\leq\frac{\nu}%
{2}||\nabla e||^{2}+\frac{1}{\nu}||\nabla u||^{2}||e||^{2}.
\]
Thus,%
\[
\frac{d}{dt}||e||^{2}+\left[  \nu-2C_{1}^{2}H^{2}\chi\right]  ||\nabla
e||^{2}+\left[  \chi-2\nu^{-1}||\nabla u||^{2}\right]  ||e||^{2}\leq0
\]
Provided the term in the first bracket is non-negative, the condition
(\ref{eqResult1}), we proceed as follows. Denote $a(t):=\chi-2\nu^{-1}||\nabla
u||^{2}\in L^{1}(0,T),A(T)=\int_{0}^{T}a(s)ds,$ $0\leq y(t)=||e||^{2}$ then
$y^{\prime}+a(t)y\leq0$. Thus, $y(T)\leq\exp\left[-A(T)\right]y(0)$. Thus, provided
(\ref{eq2dChi}), i.e.,
\[
\chi-\frac{2}{\nu}\frac{1}{T}\int_{0}^{T}||\nabla u||^{2}dt\geq\chi_{0}>0,
\]
it follows that $\exp\left[-A(T)\right]\leq\exp\left[-\chi_{0}T\right]$ and the the error
$\rightarrow0$\ uniformly in $T$ as $\chi\rightarrow\infty$ and uniformly in
$\chi$\ as $T\rightarrow\infty$.

\textbf{The 3d case:} Starting at%

\[
\frac{1}{2}\frac{d}{dt}||e||^{2}+(e\cdot\nabla u,e)+\left[  \nu-2C_{1}%
^{2}H^{2}\chi\right]  ||\nabla e||^{2}+\chi||e||^{2}=0.
\]
By the 3d Ladyzhenskaya inequality and an arithmetic-geometric mean inequality
(with exponents 4 and 4/3)%
\begin{gather*}
|(e\cdot\nabla u,e)|\leq||\nabla u||||e||_{L4}^{2}\leq||\nabla u||\left(
\frac{4}{3\sqrt{3}}||e||^{1/4}||\nabla e||^{3/4}\right)  ^{2}\\
\leq\left(  \frac{16}{27}||\nabla u||||e||^{1/2}\right)  \left(  ||\nabla
e||^{3/2}\right)  \leq\frac{\nu}{2}||\nabla e||^{2}+\left(  \frac{2048}%
{19683}\nu^{-3}||\nabla u||^{4}\right)  ||e||^{2}.
\end{gather*}
Thus,%
\[
\frac{d}{dt}||e||^{2}+\left[  \nu-2C_{1}^{2}H^{2}\chi\right]  ||\nabla
e||^{2}+2\left[  \chi-\frac{2048}{19683}\nu^{-3}||\nabla u||^{4}\right]
||e||^{2}\leq0.
\]
Following the $2d$ case, the error satisfies $||e(T)||$ $\leq\exp\left[
-A(T)\right]  ||e(0)||$ where $A(T)$ $=2\chi-\frac{2048}{19683}\nu^{-3}%
\frac{2}{T}\int_{0}^{T}||\nabla u||^{4}dt.$ This implies $||e(T)||\leq
e^{-\chi_{0}T}||e(0)||$, provided $\chi$\ is large enough, (\ref{eq3dChi}),
that%
\[
2\left[  \chi-\frac{2048}{19683}\nu^{-3}\frac{1}{T}\int_{0}^{T}||\nabla
u||^{4}dt\right]  \geq\chi_{0}>0.
\]
\end{proof}

\subsection{Interpreting the Conditions}\label{aprior-analysis} 
In the analysis, conditions on $H$ and $\chi$ arise naturally for uniform in
$T$ convergence. The question we now investigate for higher Reynolds number, turbulent flows is 
\begin{center}
\textit{How severe are these conditions for practical problems?}
\end{center}

Let the large scale length and velocity be denoted $L, U$, e.g., $L=|\Omega
|^{1/d}$ and $U^{2}=\lim\sup_{T\to \infty}\frac{1}{T}\int_{0}^{T}\frac{1}{|\Omega|}%
\int_{\Omega}|u|^{2}dxdt$. The Reynolds number is $\mathcal{R}e=LU/\nu$. The
large-scale turnover time is $T^{\ast}=L/U$. The constant $C_{1}$ in
(\ref{eqInterpError}) is dimensionless. $H$ has units of length and $\chi$\ of
$time^{-1}$.

In 2d the conditions are%
\[
\nu-2C_{1}^{2}H^{2}\chi\geq0\text{ and }\chi-\frac{1}{\nu}\frac{2}{T}\int
_{0}^{T}||\nabla u||^{2}dt\geq\chi_{0}>0.
\]
For fully developed turbulent flows in $2d$ with periodic boundary conditions,
Alexakis and Doering \cite{alexakis2006energy} have proven the following, consistent with
phenomenology,%
\begin{gather*}
\lim\sup_{T\rightarrow\infty}\frac{1}{T}\int_{0}^{T}\frac{\nu}{L^{2}}||\nabla
v||^{2}dt\leq k_{f}U^{3}\mathcal{R}e_{f}^{-1/2}\text{ where }\mathcal{R}%
e_{f}=\frac{U}{\nu k_{f}}\text{ and}\\
\text{ }k_{f}=\text{ average energy-input wave number, }\mathcal{R}e_{f}%
=\frac{U}{\nu k_{f}}=\mathcal{R}e(Lk_{f})^{-1}.
\end{gather*}
In 2d, phenomenology interprets the $\chi-$condition to mean%
\begin{gather*}
\chi>\frac{2}{\nu}\frac{1}{T}\int_{0}^{T}||\nabla u||^{2}dt\simeq\frac{2L^{2}%
}{\nu^{2}}k_{f}U^{3}\left(  \frac{U}{\nu k_{f}}\right)  ^{-1/2}\simeq
2k_{f}U(Lk_{f})^{+1/2}\mathcal{R}e^{+3/2}\\
\chi T^{\ast}\gtrsim 2(Lk_{f})^{+3/2}\mathcal{R}e^{+3/2}>>1.
\end{gather*}
The $H-$condition, $\nu-2C_{1}^{2}H^{2}\chi\geq0$ , can be rearranged to read%
\[
\left(  \frac{H}{L}\right)  ^{2}\leq\left(  2C_{1}^{2}\right)  ^{-1}%
\mathcal{R}e^{-1}\left(  \chi T^{\ast}\right)  ^{-1}.
\]
Since $\chi\gtrsim\mathcal{O}(\mathcal{R}e^{3/2})$ the $H-$condition becomes
$\left(  H/L\right)  ^{2}\leq\mathcal{O}(\mathcal{R}e^{-5/2})$. These
conditions on $H$ and $\chi$\ are severe.

In 3d, the scaling becomes worse. For fully developed turbulent flows in 3d
with periodic boundary conditions, the following has been proved by Doering and Foias \cite{doering2002energy}, again consistent with turbulent phenomenology in the large $T$ limit
\[
\frac{1}{T}\int_{0}^{T}\frac{\nu}{L^{3}}||\nabla v||^{2}dt\simeq U^{3}/L.
\]
Under the (often called a \textit{spherical cow assumption}) condition that
\[
\sqrt[4]{\frac{1}{T}\int_{0}^{T}||\nabla u||^{4}dt}\simeq\sqrt[2]{\frac{1}%
{T}\int_{0}^{T}||\nabla u||^{2}dt}%
\]
we then estimate
\[
\frac{1}{T}\int_{0}^{T}||\nabla u||^{4}dt\simeq\left(  L^{3}\nu^{-1}\frac
{1}{T}\int_{0}^{T}\frac{\nu}{L^{3}}||\nabla u||^{2}dt\right)  ^{2}%
\simeq\left(  L^{3}\nu^{-1}\frac{U^{3}}{L}\right)  ^{2}.
\]
The condition (\ref{eq3dChi}), $\chi$\ \textit{big enough}, has the
interpretation
\[
\chi T^{\ast}\gtrsim\frac{L}{U}0.1\nu^{-5}L^{6}\frac{U^{6}}{L^{2}%
}=0.1\mathcal{R}e^{5}.
\]
The $H$ \textit{small enough} condition then becomes
\[
\left(  \frac{H}{L}\right)  ^{2}\leq\left(  2C_{1}^{2}\right)  ^{-1}%
\mathcal{R}e^{-1}\left(  \chi T^{\ast}\right)  ^{-1}\lesssim\mathcal{O}%
(\mathcal{R}e^{-6}).
\]
The 3d condition $\chi T^* \gtrsim\mathcal{O}(\mathcal{R}e^{5})$ and $H/L\lesssim\mathcal{O}
(\mathcal{R}e^{-3})$ are severe restrictions.

For problems with smooth solutions (e.g., test problems constructed by the \\ method of manufactured solutions) the $H-$condition can be altered as follows.
For a vector function $w(x)$ define $\lambda_{T}(w)$ as%
\[
\lambda_{T}(w):=\frac{||w||}{||\nabla w||},d=2 \text{ or } 3.
\]
Here $\lambda_{T}$\ represents an average length-scale of the solution,
analogous to the Taylor micro-scale (which would also involve time averaging).
Then $||\nabla e||\leq\lambda_{T}^{-1}||e||$. The term that gives rise to the
$H$-condition is now estimated by%
\[
\chi||(I-I_{H})e||^{2}\leq\chi\left(  C_{1}H||\nabla e||\right)  ^{2}\leq\chi
C_{1}^{2}H^{2}\lambda_{T}^{-2}||e||^{2}.
\]
It can now be subsumed into the $\chi$-condition (in 3d here). This gives long
time estimates when $H$ is small with respect to $\lambda_{T}(e)$ (rather than $\mathcal{R}e$) from the inequality.
\[
\frac{d}{dt}||e||^{2}+\nu||\nabla e||^{2}+2\left[  \chi\left(  1-C_{1}
^{2}\left(  \frac{H}{\lambda_{T}(e)}\right)  ^{2}\right)  -\frac
{2048}{19683}\nu^{-3}||\nabla u||^{4}\right]  ||e||^{2}\leq0.
\]

\section{Self-Adaptive Parameter Selection}\label{adapt-section}
The severe constraints, in 2d
$\chi \gtrsim \mathcal{O}(\mathcal{R}e^{3/2})$ and $\chi \gtrsim \mathcal{O}
(\mathcal{R}e^{5})$\ in 3d can possibly be improved since analysis
gives a sufficient (but possibly not sharp), worst-case lower bound. It is
therefore interesting to develop methods for self-adaptive parameter selection
that respond to the local state of the flow. The condition on $H$ is linked to
the condition on $\chi$; improving the latter improves the former.
Further, changing $H$ requires changing observation locations rather than just
picking a new user-defined parameter. We thus focus on adapting $\chi$. Reformulating nudging to improve the H-condition is an open problem.

The first algorithm is not grounded in the \'{a} priori analysis other than the result that a decrease of $||e||$\ is possible and $\|I_H(e)\|\leq \|e\|$.  It uses the fact that $||I_{H}e||$ is computable in a reasonable but heuristic way as follows.

\begin{algorithm}[H]
\caption{}
Choose $\chi_{0}>0$. Set upper safety factor, Factor $\geq 1  $ and lower tolerance, $Tol <1$. Given $\chi_{n}$ such that $||I_{H}e(t_{n})||<||I_{H}%
e(t_{n-1})||$

Pick $\chi_{n+1}=\chi_{n}$ and calculate $v(t)$ for $t_{n}<t\leq t_{n+1}.$

\textbf{[Acceptable }$\chi$\textbf{\ value]} Proceed to next step with
$\chi_{n+2}=\chi_{n+1}$ if \ \ $||I_{H}e(t_{n+1})||<||I_{H}e(t_{n})||$

\textbf{[Too small }$\chi$\textbf{\ value] }If $||I_{H}e(t_{n+1})||\geq \text{ Factor }
||I_{H}e(t_{n})||$, set $\chi_{n+1}=2.0\chi_{n}$ and repeat step.

\textbf{[Too large }$\chi$\textbf{\ value] }If $||I_{H}e(t_{n+1}%
)||\leq \text{ Tol } ||I_{H}e(t_{n})||$, set $\chi_{n+2}=0.5\chi_{n+1}$ and next step.
\end{algorithm}

The second algorithm is grounded in the \'{a} priori analysis that led to conditions (\ref{eq2dChi}), (\ref{eq3dChi}). We first
develop it in 2d. Changing the parameter values at each time step means the
nudging term is no longer autonomous. Thus we set $\chi=\chi(t)$. Since we are not addressing the $H$-condition the analysis is shortened by setting $I_{H}=I$. We begin with
\[
\frac{1}{2}\frac{d}{dt}||e||^{2}+(u\cdot\nabla u-v\cdot\nabla v,e)+\nu||\nabla
e||^{2}+\chi(t)||e||^{2}=0.
\]
For the nonlinear term we rearrange and decompose on $(u\cdot\nabla u-v\cdot\nabla v,e)=(e\cdot\nabla
v,e)$ rather than $(e\cdot \nabla u,e)$ since $v$ is computed. Following the subsequent steps (making this one
change) we obtain
\[
\frac{d}{dt}||e||^{2}+\nu||\nabla e||^{2}+\left[  2\chi(t)-\nu^{-1}||\nabla
v||^{2}\right]  ||e||^{2}\leq0.
\]
Denote $a(t):=2\chi(t)-\nu^{-1}||\nabla v||^{2},0\leq y(t)=||e||^{2} $
then $y^{\prime}+a(t)y\leq0$. Thus, the error at time $T$ satisfies%
\[
||e(T)||\leq\exp\left[  -\left(  \frac{1}{T}\int_{0}^{T}\chi(t)dt-\frac
{1}{2\nu}\frac{1}{T}\int_{0}^{T}||\nabla v||^{2}dt\right)  T\right]  ||e(0)||
\]
and the contribution to the evolution of the error of one-time step of size
$\tau$ satisfies%
\[
||e(t+\tau)||\leq\exp\left[  -\left(  \frac{1}{\tau}\int_{t}^{t+\tau}%
\chi(s)ds-\frac{1}{2\nu}\frac{1}{\tau}\int_{t}^{t+\tau}||\nabla v(s)||^{2}%
ds\right)  \tau\right]  ||e(t)||
\]

This suggests the following algorithm with piecewise constant $\chi$.

\begin{algorithm}[H]
\caption{}{}
Choose $\chi_{0}>0$. Given $\chi_{n}$ such that
\[
\chi_{n}-\frac{1}{2\nu}\frac{1}{\tau}\int_{t_{n-1}}^{t_{n}}||\nabla
v(s)||^{2}ds\geq\chi_{0}>0.
\]
Pick $\chi_{n+1}=\chi_{n}$ and calculate $v(t)$ for $t_{n}<t\leq t_{n+1}.$

\textbf{[Acceptable }$\chi$\textbf{\ value]} Proceed to next step with
$\chi_{n+2}=\chi_{n+1}$ if
\[
2\chi_{0}\geq\chi_{n+1}-\frac{1}{2\nu}\frac{1}{\tau}\int_{t_{n}}^{t_{n+1}%
}||\nabla v(s)||^{2}ds\geq\chi_{0}>0.
\]

\textbf{[Too large }$\chi$\textbf{\ value] }If
\[
\chi_{n+1}-\frac{1}{2\nu}\frac{1}{\tau}\int_{t_{n}}^{t_{n+1}}||\nabla
v(s)||^{2}ds>2\chi_{0}
\]
Proceed to the next step after resetting
\[
\chi_{n+2}=1.1\chi_{0}+\frac{1}{2\nu}\frac{1}{\tau}\int_{t_{n}}^{t_{n+1}%
}||\nabla v(s)||^{2}ds.
\]

\textbf{[Too small }$\chi$\textbf{\ value] }If%
\[
\chi_{0}>\chi_{n+1}-\frac{1}{2\nu}\frac{1}{\tau}\int_{t_{n}}^{t_{n+1}}||\nabla
v(s)||^{2}ds
\]
Repeat step after resetting%
\[
\chi_{n+1}=1.1\chi_{0}+\frac{1}{2\nu}\frac{1}{\tau}\int_{t_{n}}^{t_{n+1}%
}||\nabla v(s)||^{2}ds.
\]

\end{algorithm}

In 3d the only change is to replace
\begin{equation*}
\frac{1}{2\nu}\frac{1}{\tau}\int_{t_{n}%
}^{t_{n+1}}||\nabla v(s)||^{2}ds \text{ by } 
\frac{2048}{19683}\nu^{-3}\frac{1}{\tau
}\int_{t_{n}}^{t_{n+1}}||\nabla u||^{4}dt.
\end{equation*}
In both cases, the integrals over
$t_{n}<t<t_{n+1}$ are evaluated by quadrature using the trapezoidal
approximation and the approximate velocities at $t_{n},t_{n+1}$.

\section{Numerical Tests}\label{sec: tests}
In this section, we conduct three tests to validate the adaptive algorithms. In the first test, we calculate the rate of convergence with a test with exact solutions. Then we test the adaptive algorithms with a complex flow between offset cylinders from \cite{jiang2016algorithms} with a higher Reynolds number. In the third test, we use the test case from \cite{leotest}, which involves flows past a flat plate with a $Re$ of $50$.

We note that the two adaptive strategies have similar but different aims. The first aims to enforce $\|I_H  (e(t+\tau))\|\leq\|I_H (e(t))\|$. The second aims to enforce $\|e(t+\tau)\|\leq \|e(t)\|$. Since $\|I_H(e)\|\leq \|e\|$ we expect the first to yield smaller $\chi$-values than the second. 
\subsection{Test of Accuracy}
\begin{figure}[H]
    \centering
\includegraphics[width=0.75\linewidth]{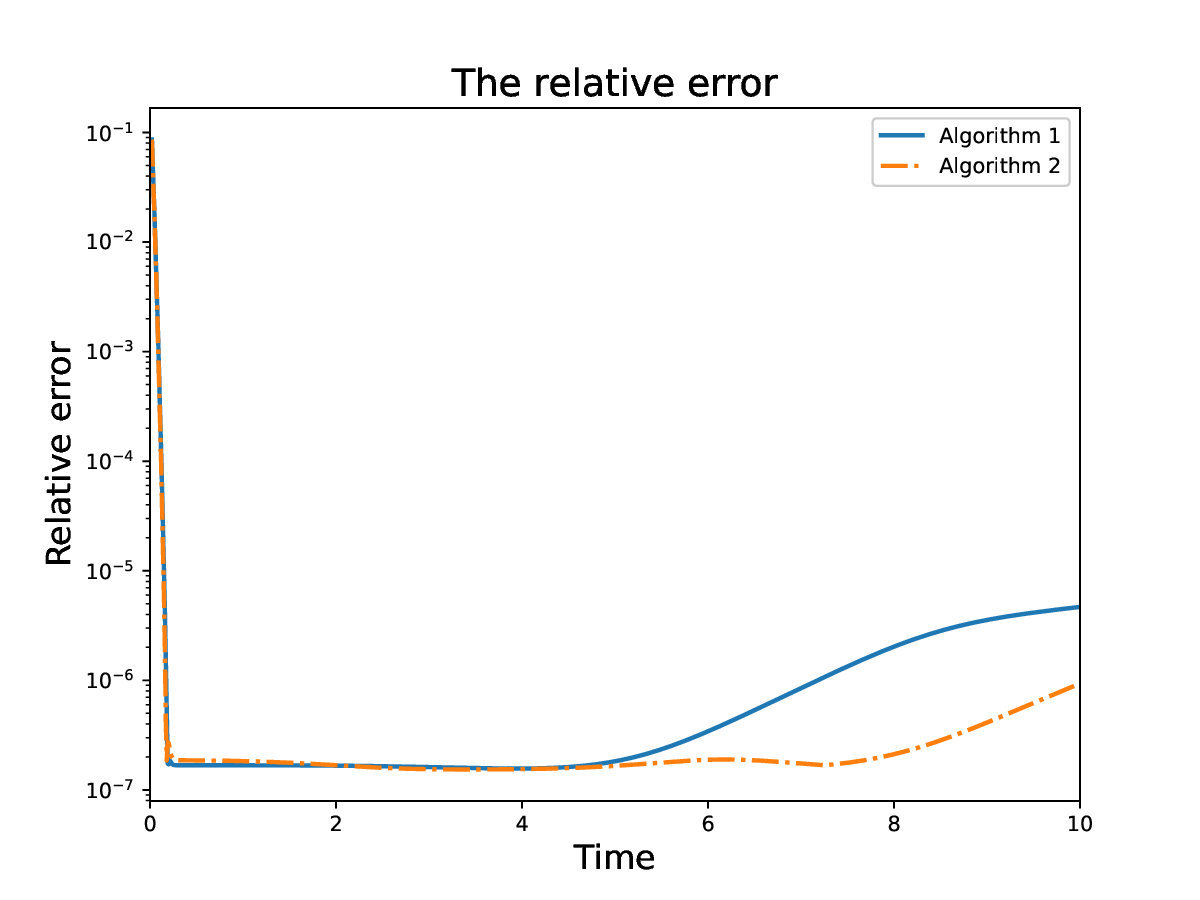}
    \caption{Relative errors. In a longer time, Algorithm 2 gives smaller errors.}
    \label{fig: exact-error}
\end{figure}

\begin{figure}[H]
    \centering
\includegraphics[width=0.75\linewidth]{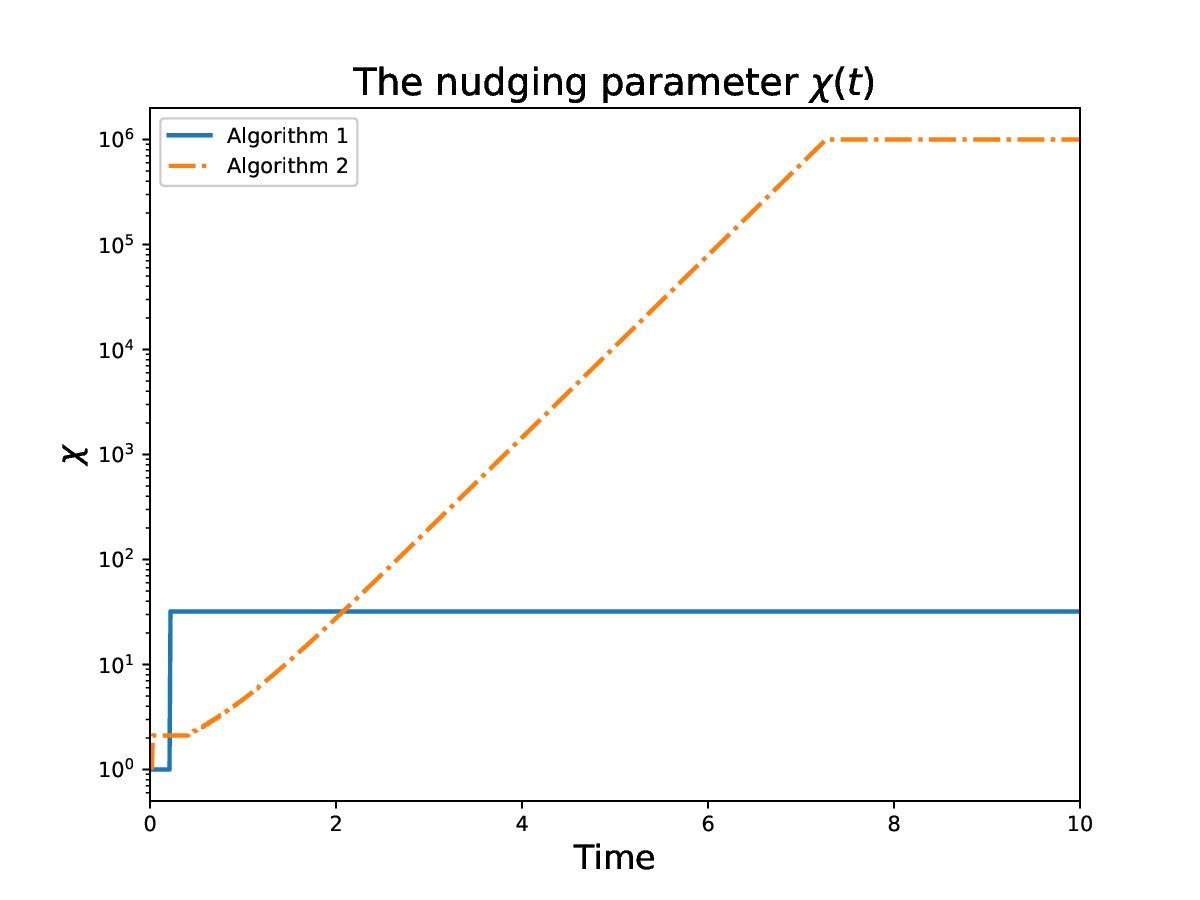}
    \caption{$\chi$ in long time. $\chi$ values in Algorithm 2 reach the maximum value.}
    \label{fig: exact-chi}
\end{figure}

In the numerical example, we verify the theoretical temporal convergence results with the BDF2 (second-order backward differentiation formula) time discretization described and analyzed in Section 3 of Larios, Rebholz, and Zerfas, \cite{larios2019global}. The following exact solution for the NSE is considered in the domain $\Omega=(0,1)\times(0,1)$. The velocity and the pressure are as follows:
\begin{equation*}
    u(x,y,t)=e^t(\cos y,\sin x), \ \text{ and }
    p(x,y,t)=(x-y)(1+t).
\end{equation*}
We insert these in the NSE to calculate the body force $f(x,t)$. We choose $Re=1$. We create a fine mesh resolution that consists of $55554$ degrees of freedom (dofs) to isolate the convergence rates for $\Delta t $ of the errors by making the spacial errors small. We run the code in the time interval $[0,2]$. We use a Scott-Vogelius finite element pair with a barycenter refined mesh. We observed similar results with Taylor-Hood elements and a skew-symmetrized nonlinearity. The measurements were made by $L^2(\Omega)$ norm. In Algorithm 1, we set a safety factor and a lower Tolerance as follows. The factor is user-defined to avoid $\chi$ growing too big when the projection error grows slightly. 
\begin{equation*}
\text{Doubling when } \|I_H(e(t+\tau))\|> \text{ Factor } \|I_H(e(t))\|.
\end{equation*}
We set Factor $= 1.3$.
\begin{equation*}
\text{Halving when } \|I_H(e(t+\tau))\|< \text{ Tol} \|I_H(e(t))\|,
\end{equation*}
where Tol is also defined by the user and Tol=$0.2$ in this test. 
\begin{table}[H]
    \centering
    \begin{subtable}{0.4\linewidth} %
     \centering
        \begin{tabular}{||c|c|c|c||}
            \hline
            $\Delta t$ & $\|u-v\|$ & rate & $\chi_{\max}$ \\
            \hline
            1       & 0.0052   & -    & 1    \\
            1/2     & 0.0013   & 2.01 & 1    \\
            1/4     & 0.00026  & 2.32 & 32   \\
            1/8     & 5.4e-5   & 2.26 & 64   \\
            1/16    & 1.4e-5   & 1.94 & 64   \\
            1/32    & 1.8e-6   & 2.92 & 256  \\
            \hline
        \end{tabular}
        \caption{Algorithm 1}
        \label{tab:accu2}
    \end{subtable}
    \hfill
    \begin{subtable}{0.4\linewidth}
 \begin{tabular}{||c|c|c|c||}
            \hline
            $\Delta t$ & $\|u-v\|$ & rate & $\chi_{\max}$ \\
            \hline
            1       & 0.0042   & -    & 16  \\
            1/2     & 0.0010   & 2.07 & 19   \\
            1/4     & 0.00028  & 1.83 & 23   \\
            1/8     & 7.8e-5  & 1.84 & 20   \\
            1/16    & 2.0e-5   & 1.96 & 23   \\
            1/32    & 5.2e-6   & 1.94 & 27   \\
            \hline
        \end{tabular}
        \caption{Algorithm 2}
        \label{tab:accu1}
       
    \end{subtable}
    \caption{Errors, rates of convergence, and values of the maximum value of $\chi$ through adaptation for Algorithm $1$ and $2$. We obtain the second-order error rates in time.}
    \label{tab:error-main}
\end{table}
In Table \ref{tab:error-main}, we observe second-order convergence in $\Delta t$ for velocity. This is an optimal rate and is evidence that the implementation is correct. We also present the values of the maximum value of $\chi$ through adaptation to understand how large the generated $\chi$ values become. The small $\chi$ values seem to be in response to the solution's smaller $\lambda_T$. According to the results of this accuracy test, theoretical expectations are compatible with numerical results, except that, oddly, the $\chi$-values of Algorithm $1$ are larger than Algorithm $2$. This effect reverses over a longer time.   
Due to this reversal of the expected $\chi$ values, we monitor the behavior of the relative errors and $\chi$ values for a longer time, final time $T=10$, see Figure \ref{fig: exact-error} and Figure \ref{fig: exact-chi}. The relative errors are $\mathcal{O}(10^{-6})$ for a long time. Then they start growing, but the growth speed is slowing down and the errors appear to be saturating around $10^{-6}$ similar to logistic growth. Logistic growth is a long-accepted model of error growth in fluid dynamics, e.g., Lorenz \cite{lorenz1985growth}. For Algorithm 1, the adapted $\chi$ value is less than $50$, and for Algorithm $2$, $\chi$ is growing until it reaches the maximum value $10^{6}$. Algorithm $1$ provided smaller $\chi$-values is consistent with expectations, based on $\|I_H(e)\|\leq \|e\|$.

\subsection{Flow between offset cylinders}
In this example, we are testing the adaptive nudging algorithms with a higher $Re$ flow with complex features and without exact solutions. The domain is a disk with a smaller off-center obstacle inside. Let the outer circle radius $r_1 = 1$ and the inner circle radius $r_2 = 0.5$, with centers $c = (c_1, c_2)= (\frac{1}{2},0)$. The domain is
\begin{equation*}
\Omega = \{ (x,y): x^2 + y^2 \leq r_1^2 \text{ and } (x-c_1)^2 + (y-c_2)^2 \geq r_2^2\}.
\end{equation*}
A rotational force drives the flow,
\begin{equation*}
f(x,y,t) = (-4 y \min(1,t) \ (1-x^2-y^2), 4x\min(1,t) \ (1-x^2-y^2))^\top,
\end{equation*}
with no-slip boundary conditions. The outer circle remains
stationary. The Delaunay algorithm generates the mesh with 75 mesh points on the outer circle and 60 mesh points on the inner circle. The final time $T = 10$, time step size $\Delta t = 0.01$, $\nu = 10^{-3},$ $L=1, U=1$, and $Re = \frac{UL}{\nu}$. Initial condition $u(x,y,0) = 0$ and the Dirichlet boundary condition $u=0$ on $\partial \Omega$. 

We use BDF2  for the time discretization and Taylor-Hood ($P2-P1$) for the velocity and pressure spaces. We apply the second-order linear extrapolation to the explicitly skew symmetrized nonlinear term $b^*(v^*, v^{n+1}, v^h)$, where
\begin{equation*}
b^*(u,v,w) = \frac{1}{2} (u\cdot \nabla v, w) -\frac{1}{2} (u\cdot \nabla w, v)\text{ and } v^*=2v^n-v^{n-1}.
\end{equation*}
In this test, we approximate the true solution $u_{obs}$ for nudging using two methods: Direct Numerical Simulation (DNS) and an Unsteady Reynolds-Averaged Navier-Stokes (URANS) model of turbulence from \cite{fang20231, han2024numerical}.

We compare the adaptive Algorithm 1, Algorithm 2, and constant $\chi = 10^4$. The initial nudging parameter for the adaptive algorithms is $\chi(0)=1$. The upper bound for $\chi$ is $\chi_{\max}=10^{6}$.  In Algorithm 1, we set $\text{Factor} = 1.1$ and $Tol=0.3$.

\subsubsection{DNS}
\begin{figure}[H]
    \centering
    \includegraphics[width=0.75\linewidth]{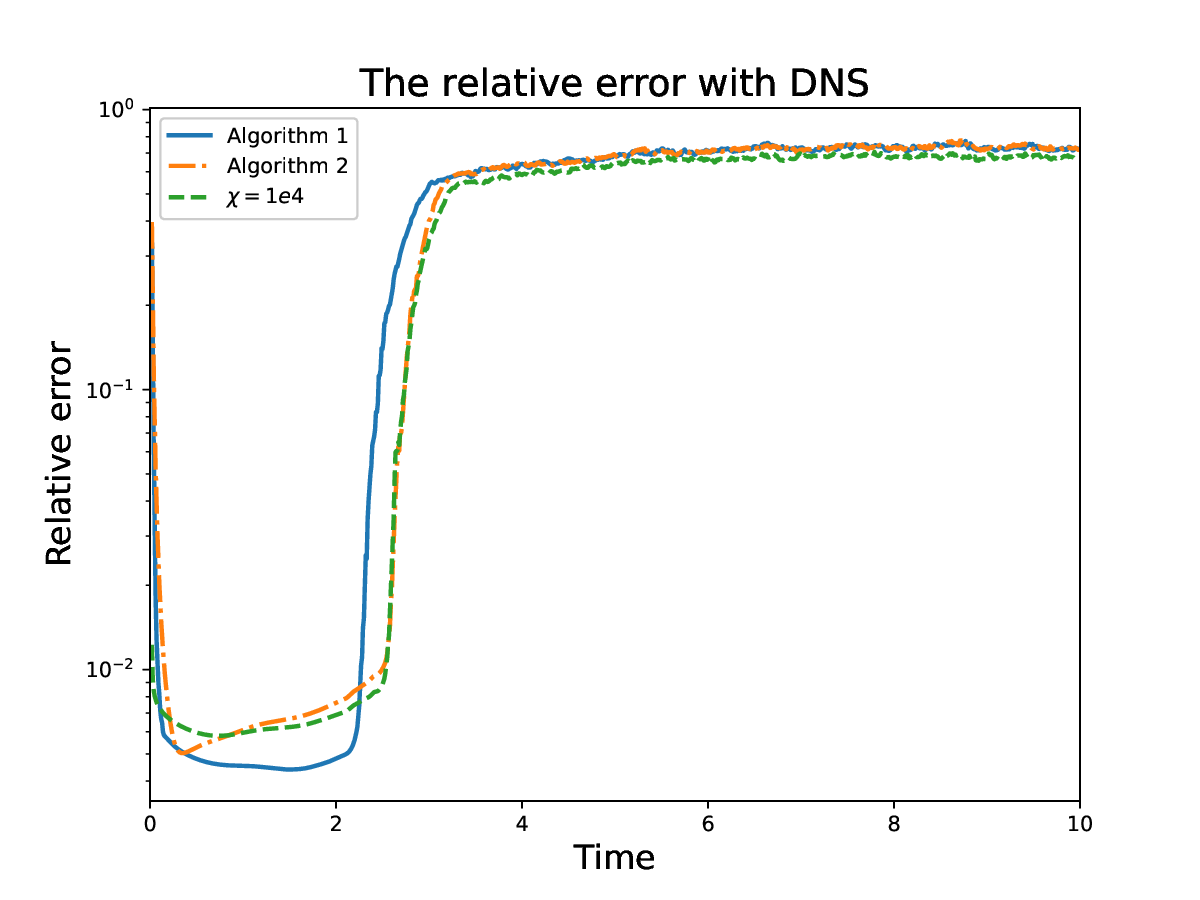}
     \caption{The relative error with the approximate true solution from DNS. The adaptive algorithms are effective for a shorter time. When the flow becomes more complex, the errors grow big and saturate at $\mathcal{O}(1)$.}
    \label{fig: DNS-error}
\end{figure}

\begin{figure}[H]
    \centering
\includegraphics[width=0.75\linewidth]{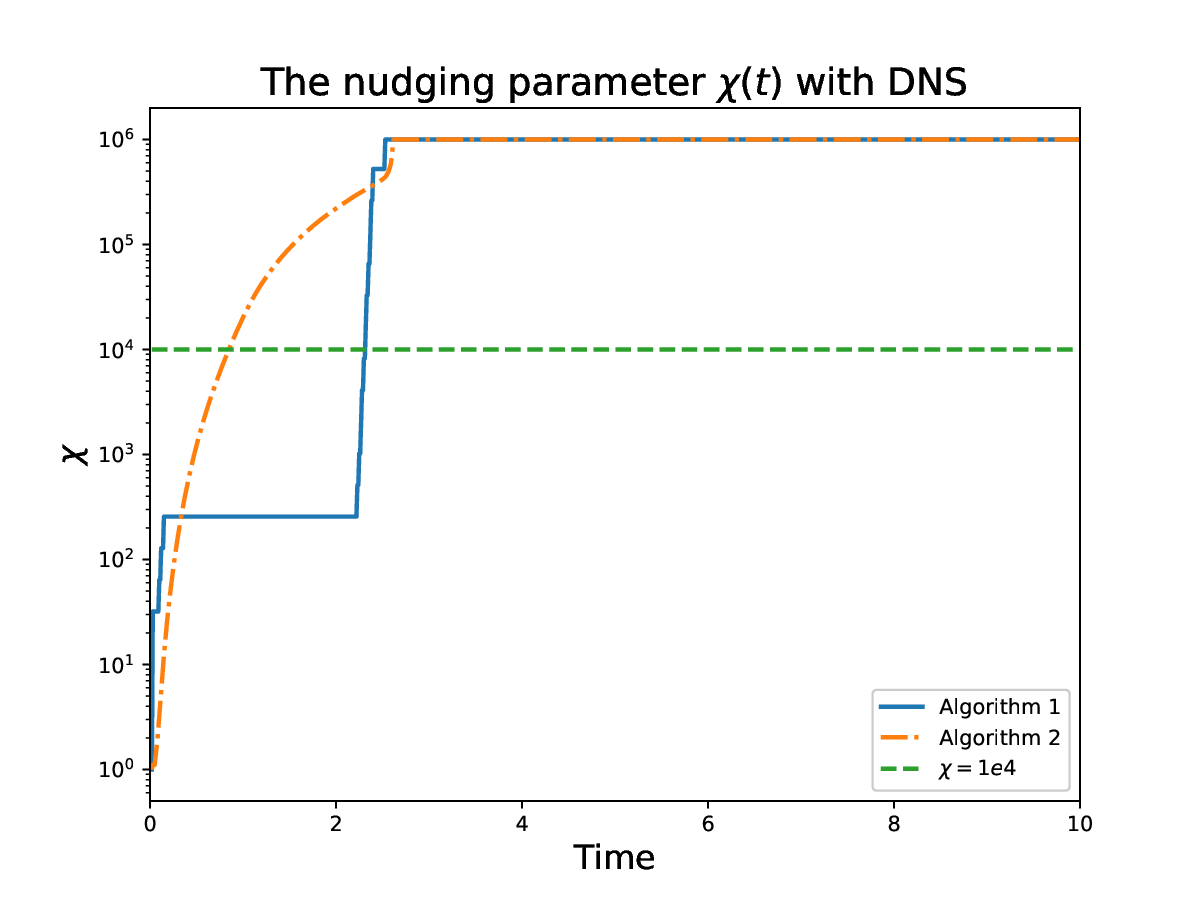}
        \caption{The plot of the nudging parameter $\chi$ values in time. Both Algorithms 1 and 2 reach the maximum $\chi$ value in a longer time.}
    \label{fig: DNS-chi}
\end{figure}

The well-resolved NSE solution is calculated on a finer mesh with $120$ mesh points on the outer circle and $96$ mesh points on the inner circle. The relative error is defined as $\frac{\|u_{DNS}-v\|}{\|u_{DNS}\|}$. In Figure \ref{fig: DNS-error}, we can see the adaptive algorithms and the constant nudging effectively decrease the relative error from $t=0$ to $t=2$. After $t=2$, the relative error grows until it saturates at $\mathcal{O}(1)$. The behavior indicates that an H-condition like (\ref{eqResult1}) for long-time convergence is necessary and not satisfied here. In Figure \ref{fig: DNS-chi}, the $\chi$ values of Algorithm $1$ and Algorithm 2 grow with time until they reach the maximum value. We observe that the $\chi$ of Algorithm 1 grows fast at first then the growth rate gets smaller. Algorithm $1$ gives a slightly better approximation from $t=0$ to $t=2$. $\chi$ for Algorithm $2$ grows at first then barely changes, then grows big due to the growth of more complex flow structures and in response to the (we conjecture) non-satisfaction of the $H$-condition. 

\subsubsection{URANS}
\begin{figure}[H]
    \centering
    \includegraphics[width=0.75\linewidth]{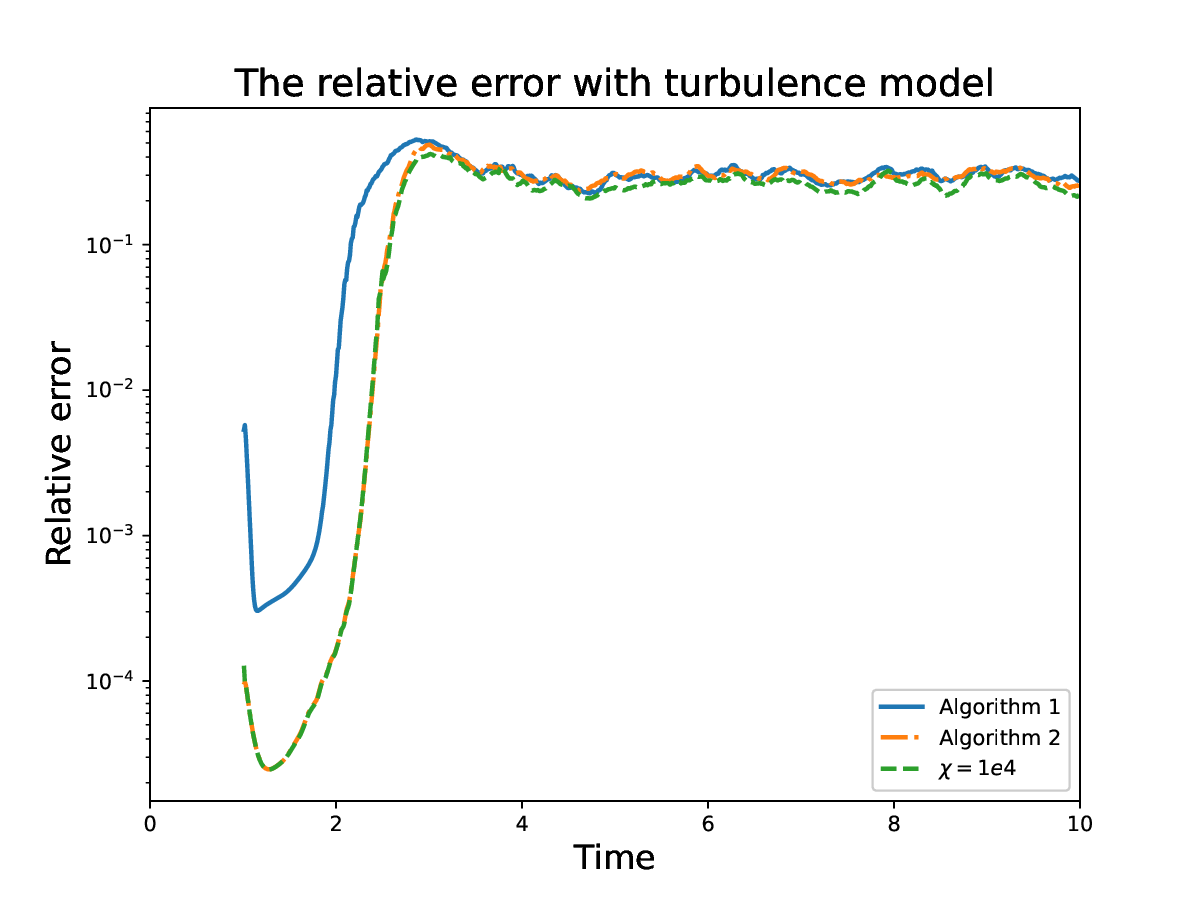}
    \caption{The relative error with the approximate true solution from a URANS model. When the flow becomes more complex, the errors grow big and saturate at $\mathcal{O}(1)$}
    \label{fig: URANS-error}
\end{figure}

\begin{figure}[H]
    \centering
    \includegraphics[width=0.75\linewidth]{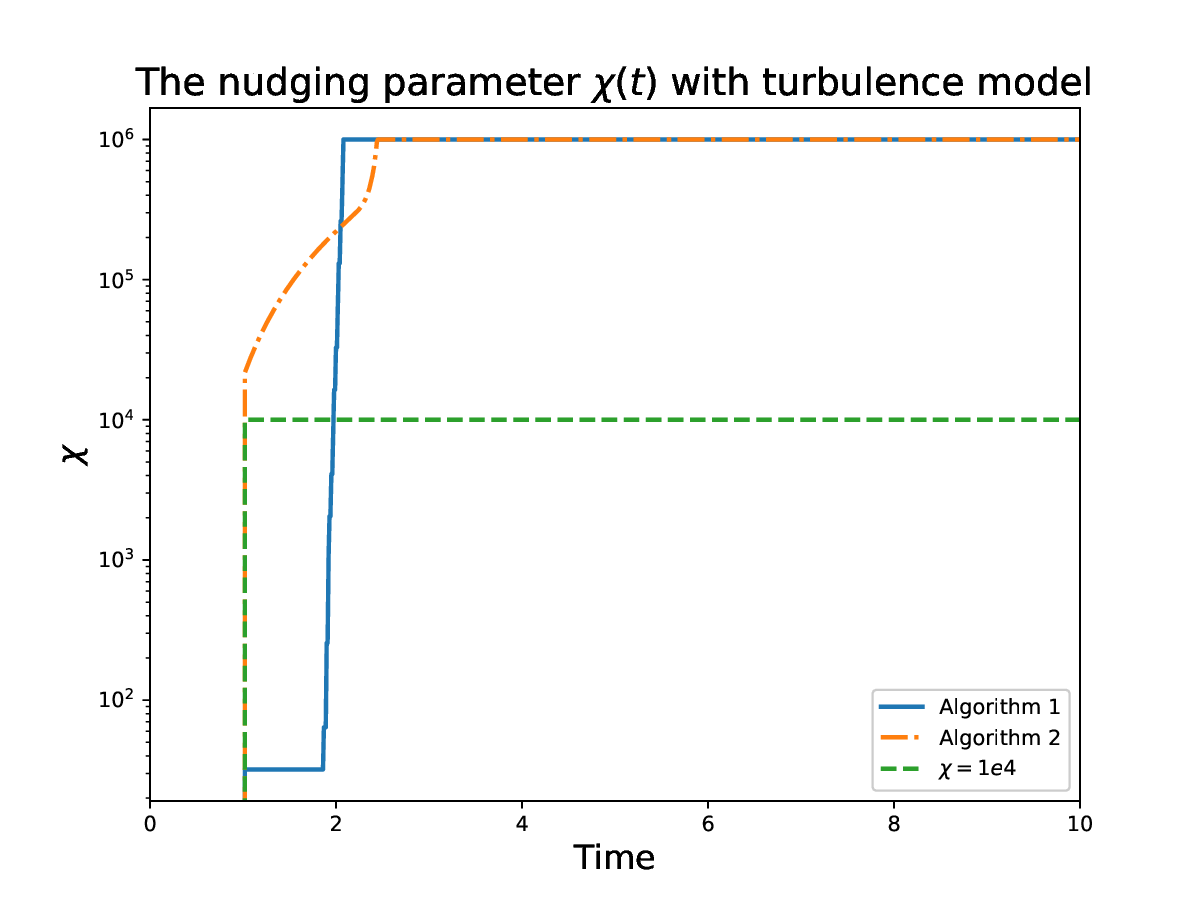}
    \caption{The plot of the nudging parameter $\chi$ values in time. Both Algorithms 1 and 2 reach the maximum $\chi$ value by time = 2.5.}
    \label{fig: URANS-chi}
\end{figure}

We use a 1/2-equation model of URANS from \cite{fang20231, han2024numerical} to approximate the exact velocity $u$. We denote it as $u_{URANS}$. The 1/2-equation model produces velocity statistics comparable to the full 1-equation model, \cite{fang20231}. The 1/2-equation URANS model is
\begin{equation}
\begin{gathered}
\frac{\partial u}{\partial t} + u \cdot \nabla u -\nabla \cdot ([2\nu + \nu_T]\nabla^s u) +\nabla q =\frac{1}{\tau} \int_{t-\tau}^t f(x,t')\, dt' \text{ and } \nabla \cdot u =0,\\
\frac{\partial k(t)}{\partial t} + \frac{\sqrt{2}}{2} \tau^{-1} k(t) = \frac{1}{|\Omega|} \int_{\Omega} \nu_t |\nabla^s v |^2\, dx, 
\end{gathered}
\end{equation}
where $\nu_t = \sqrt{2}\mu k(t) (\kappa \frac{y}{L})^2 \tau$ for no-slip and shear boundaries. Here, $y$ is the wall-normal distance. We choose the calibration constant $\mu=0.55$, the time window $\tau=0.1$, and the von Karman constant $\kappa = 0.41$. We turn on the k-equation at $t^\star=1$. We initialize the 1/2-equation model with values from Wilcox \cite{wilcox1998turbulence}, also see Kean, Layton, and Schneier \cite{kean2022prandtl}:
\begin{equation*}
k(1)=\frac{1}{|\Omega|} 1.5 \ I^2 \int_{\Omega} |v(x,1)|^2\,dx , \text{ and }
I =\text{ turbulence intensity }\approx 0.16 Re^{-1/8}. 
\end{equation*} 
The relative error is defined as $\frac{\|u_{URANS}-v\|}{\|u_{URANS}\|}$ starting at $t^*=1$. In Figure \ref{fig: URANS-error}, we can see the adaptive algorithms and the constant nudging effectively decrease the relative error from $t=1$ to $t=2$. After $t=2$, the relative error grows fast until it reaches $\mathcal{O}(1)$ as in the previous test, again suggesting $H$-condition violation. In Figure \ref{fig: URANS-chi}, the $\chi$ values of Algorithm $1$ and Algorithm 2 grow with time until they reach the maximum value. We observe that the $\chi$ of Algorithm 1 grows fast at first, then the growth rate gets smaller, and $\chi$ for Algorithm $2$ at the beginning barely changes and sharply grows after $t=2$. These two tests suggest that both a $\chi$-condition and an H-condition are essential for long-time convergence. This test also suggests that data from a turbulence model does not alleviate non-satisfaction of nudging's $H$-condition.

\subsection{Flow over a flat obstacle} 
\begin{figure}[H]
     \centering
        \centering
    \includegraphics[width=0.8\linewidth]{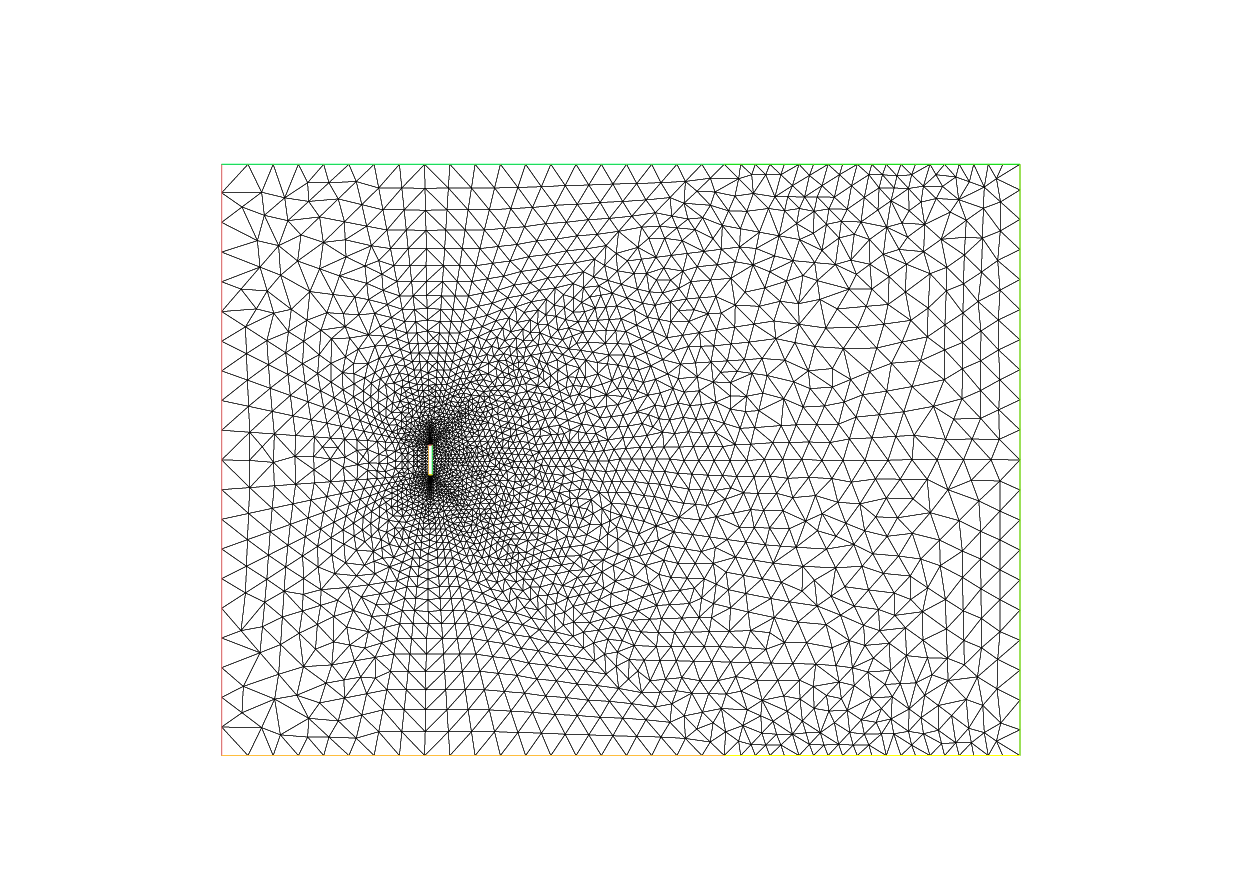}
         \caption{Finer mesh. The total dof for fine mesh is 27373 and for the coarse mesh is 15037.}
         \label{fig: Leo-mesh-finer}
\end{figure}

We adopt this test from \cite{leotest}. The computational domain is a $[-7,20]\times [0,20]$ rectangular channel with $0.125 \times 1$ flat plate obstacle so that $|\Omega|\approx 540$. The statistics are relative to solution size which also makes them area independent. The inflow velocity is set with $u_{in}= <1, 0>^T$ and
no force is applied. No-slip boundary conditions are applied on walls and the plate whereas a weak zero-traction boundary condition is enforced on the outflow boundary as in \cite{leotest}.
 We choose $\Delta t= 0.02$ and the final time $T=81$, three times turnover time. The flow is at rest at $t=0$. We approximate the exact velocity $u$ via finer mesh as seen in Figure \ref{fig: Leo-mesh-finer} and solve the $v$ equation with coarser mesh. The total dof for the fine mesh is $27373$ and for coarse mesh is $15037$.

\begin{figure}[H]
         \centering    \includegraphics[width=0.8\linewidth]{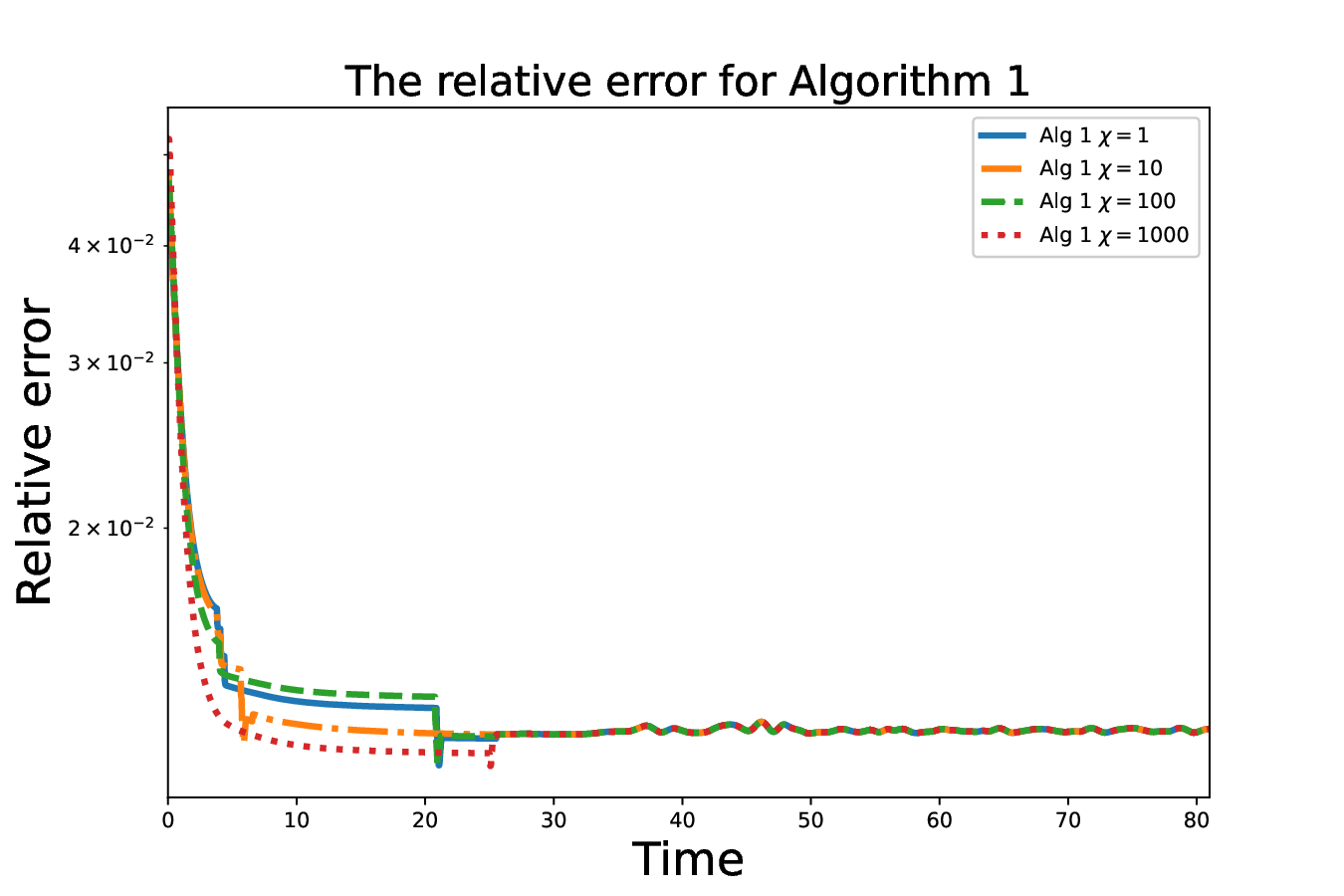}
         \caption{Algorithm 1: the relative error is dependent on the initial $\chi$ values in the first turnover time. After $\chi$ reaches the maximum value, the relative error is not sensitive to the initial $\chi$ values.}
         \label{fig: Leo-algo1-rela-err}
    \end{figure}

\begin{figure}[H]
     \centering
    \includegraphics[width=0.8\linewidth]{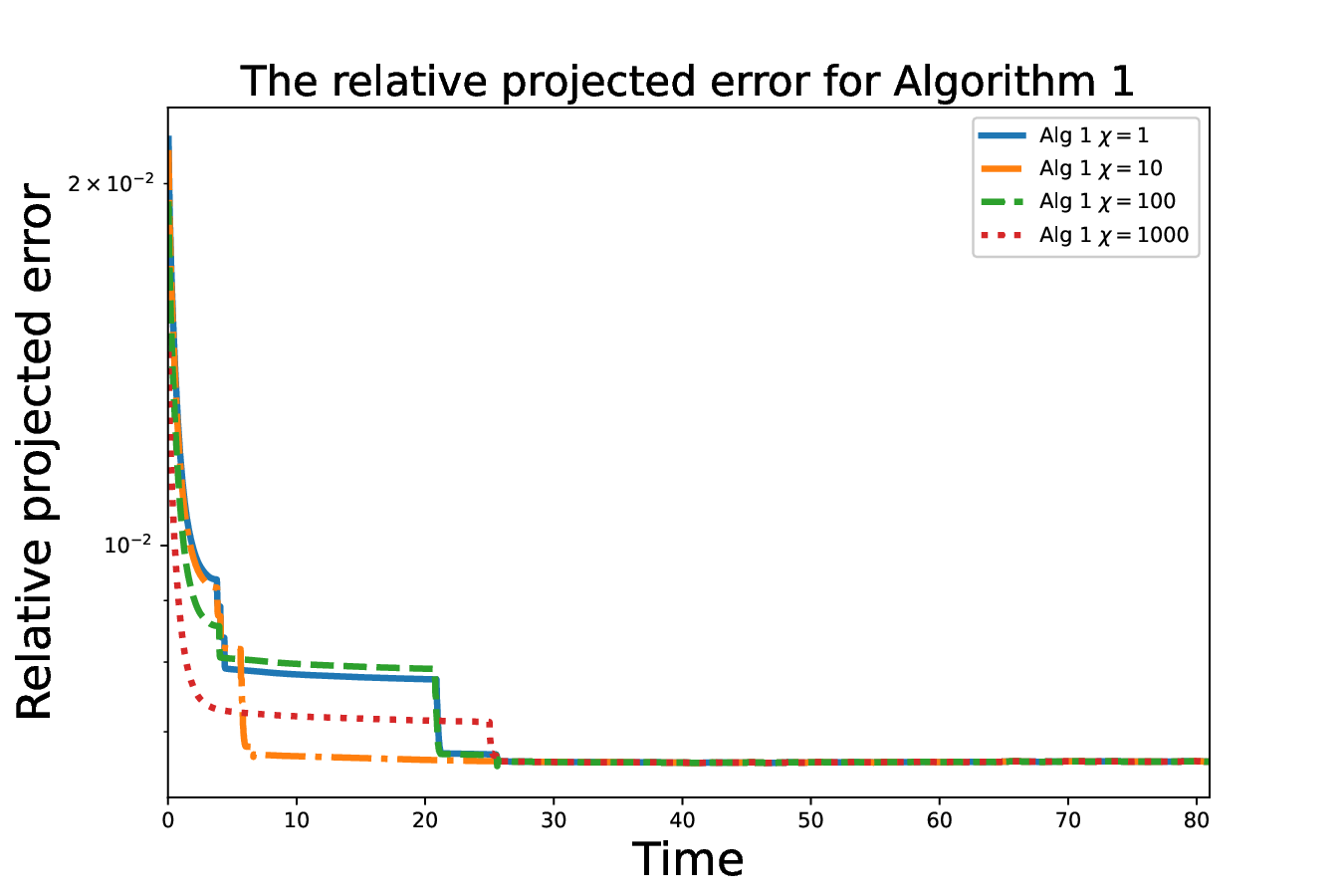}
    \caption{Algorithm 1: the projection error is dependent on the initial $\chi$ values in the first turnover time. The projection errors are monotone decreasing since the algorithm enforces $\|I_H(e(t+\tau))\|\leq \|I_H(e(t))\|$. After $\chi$ reaches the maximum value, the relative error is not sensitive to the initial $\chi$ values. }
    \label{fig: Leo-algo1-proj-err}
\end{figure}

\begin{figure}[H]
     \centering
    \includegraphics[width=0.8\linewidth]{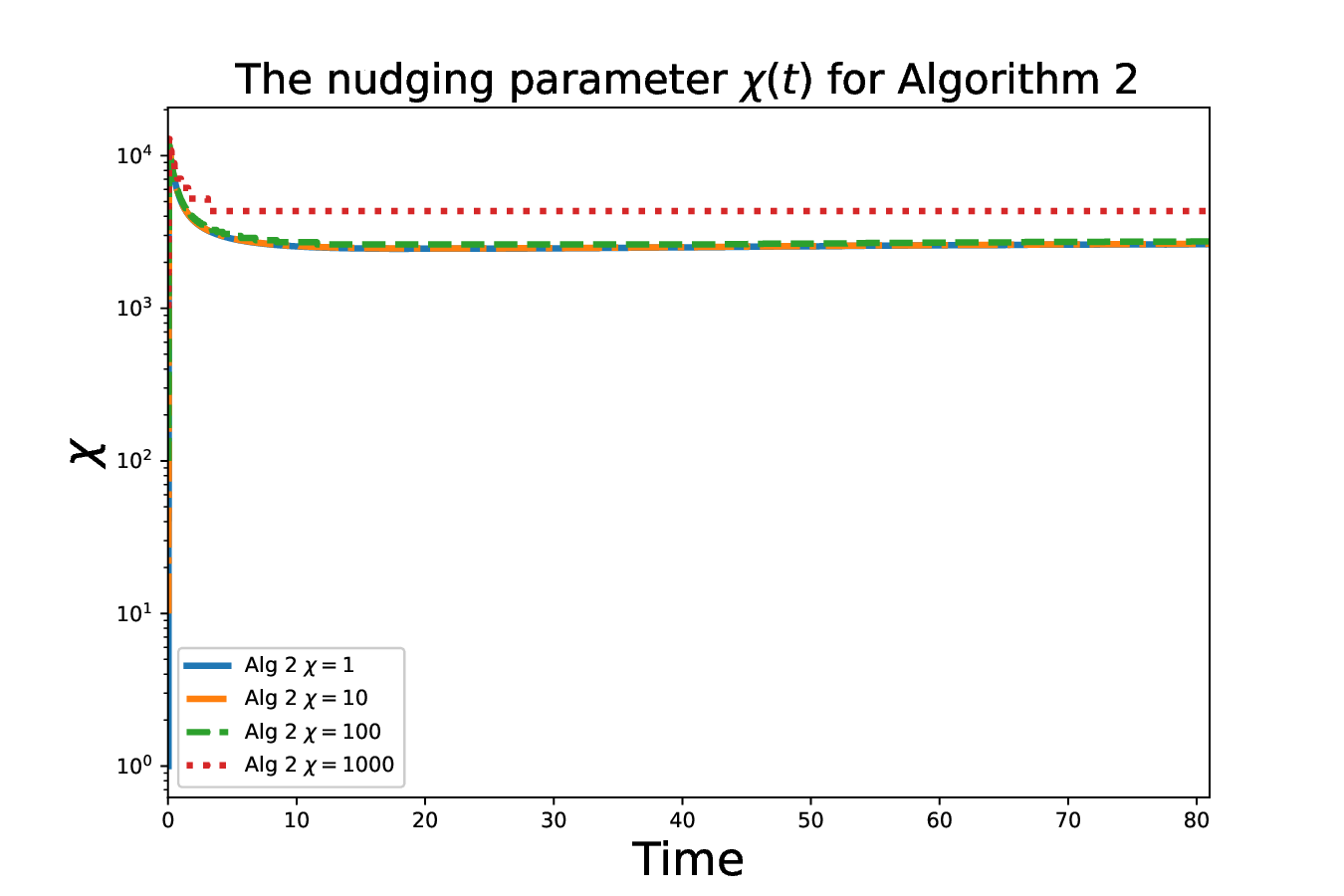}
    \caption{Evolution of $\chi$ values for Algorithm 2. $\chi$ values are adjusted quickly, and reach $10^4$ then drop out a bit and reach a steady state. For initial $\chi=1000$, the $\chi(t)$ is higher than the others.}
        \label{fig: Leo-chi-algo2}
\end{figure}

\begin{figure}[H]
         \centering    \includegraphics[width=0.8\linewidth]{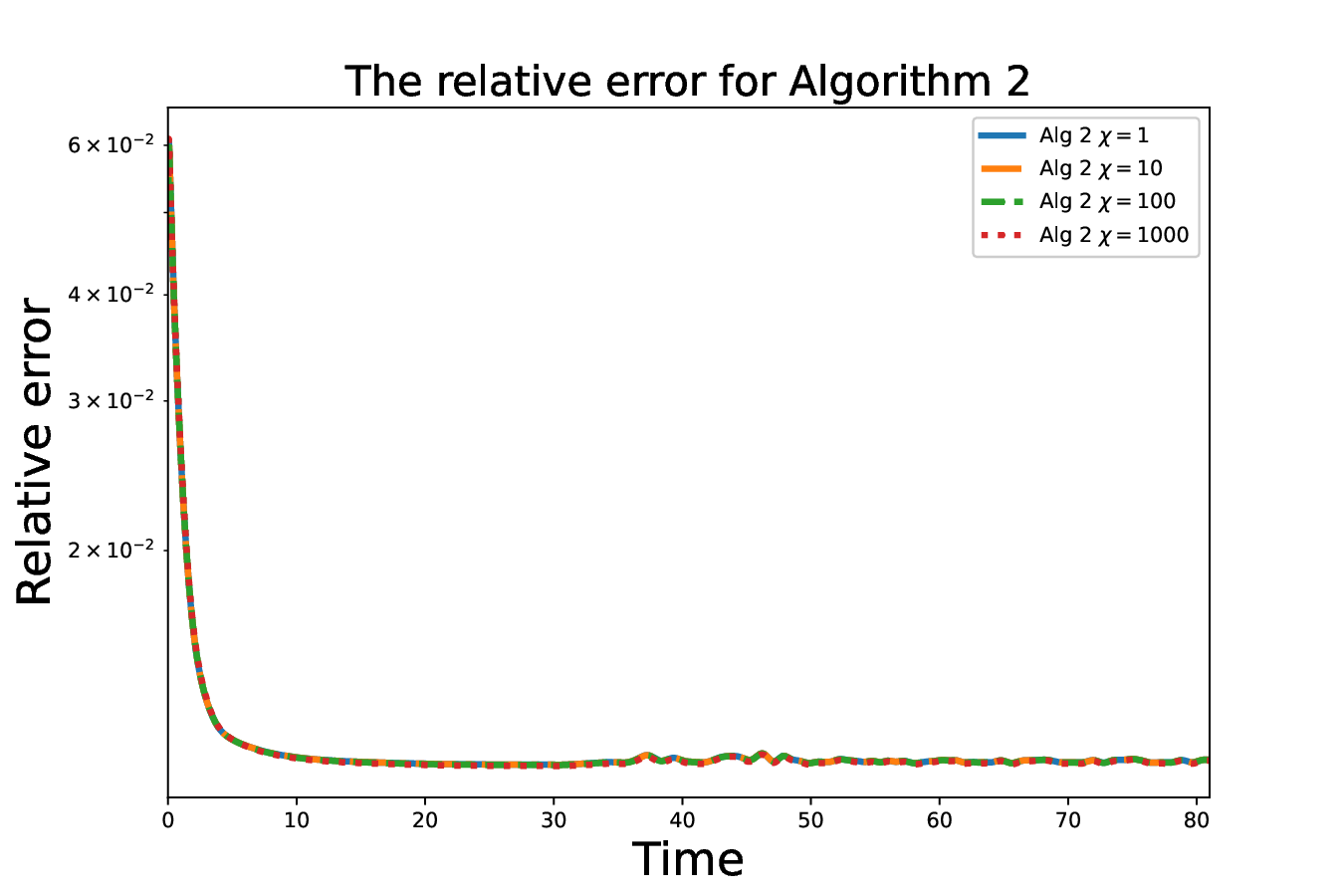}
         \caption{Algorithm 2: the relative errors are insensitive to initial $\chi$. It adapts $\chi$ value effectively to obtain smaller errors.}
         \label{fig: Leo-algo2-rela-err}
\end{figure}

\begin{figure}[H]
     \centering
    \includegraphics[width=0.8\linewidth]{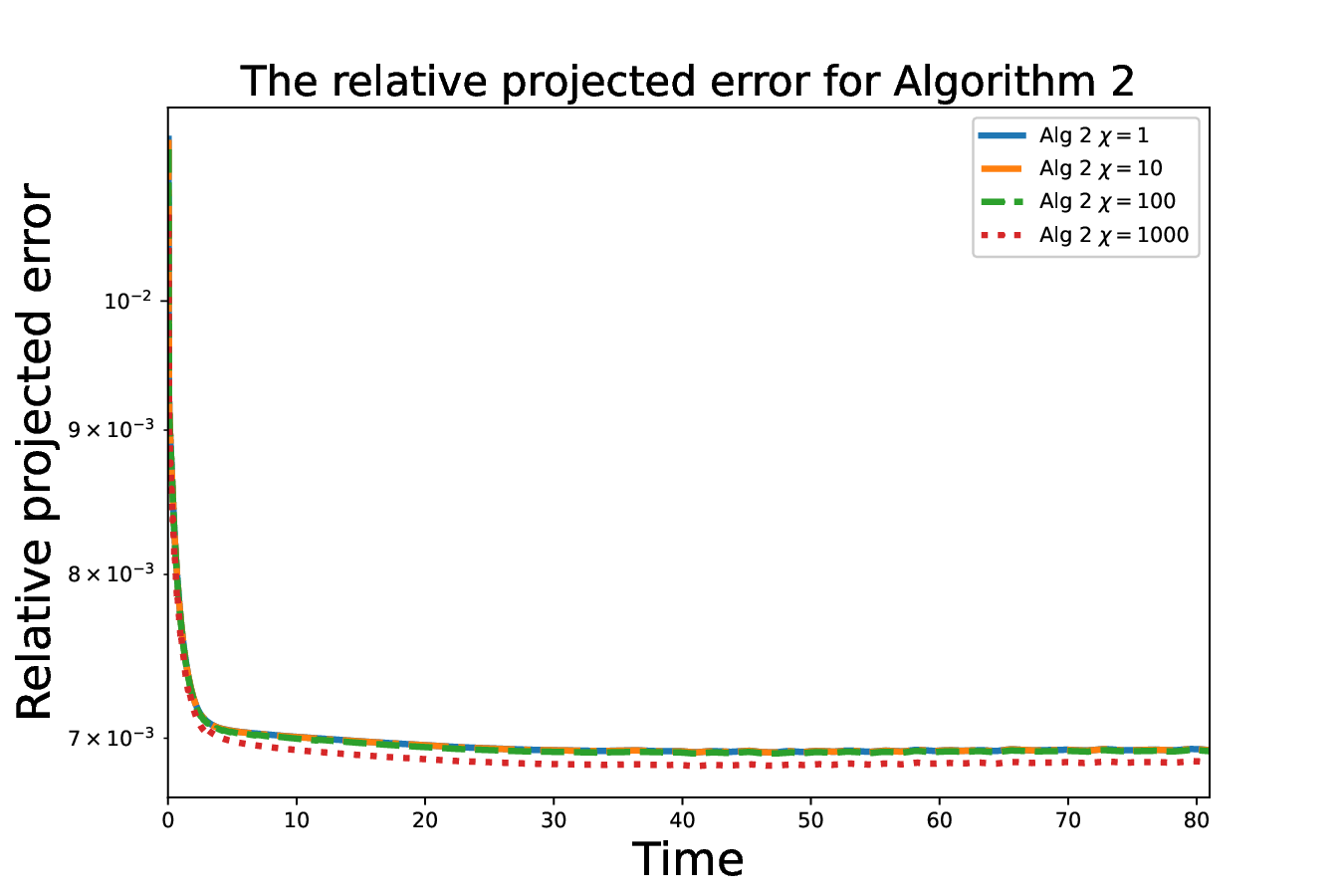}
         \caption{Algorithm 2: the projection errors are insensitive to the initial $\chi$ values. With a higher $\chi$ value for initial $\chi=1000$, it has a slightly smaller projection error.}
         \label{fig: Leo-algo2-proj-err}
\end{figure}

\begin{figure}[H]
         \centering
    \includegraphics[width=0.8\linewidth]{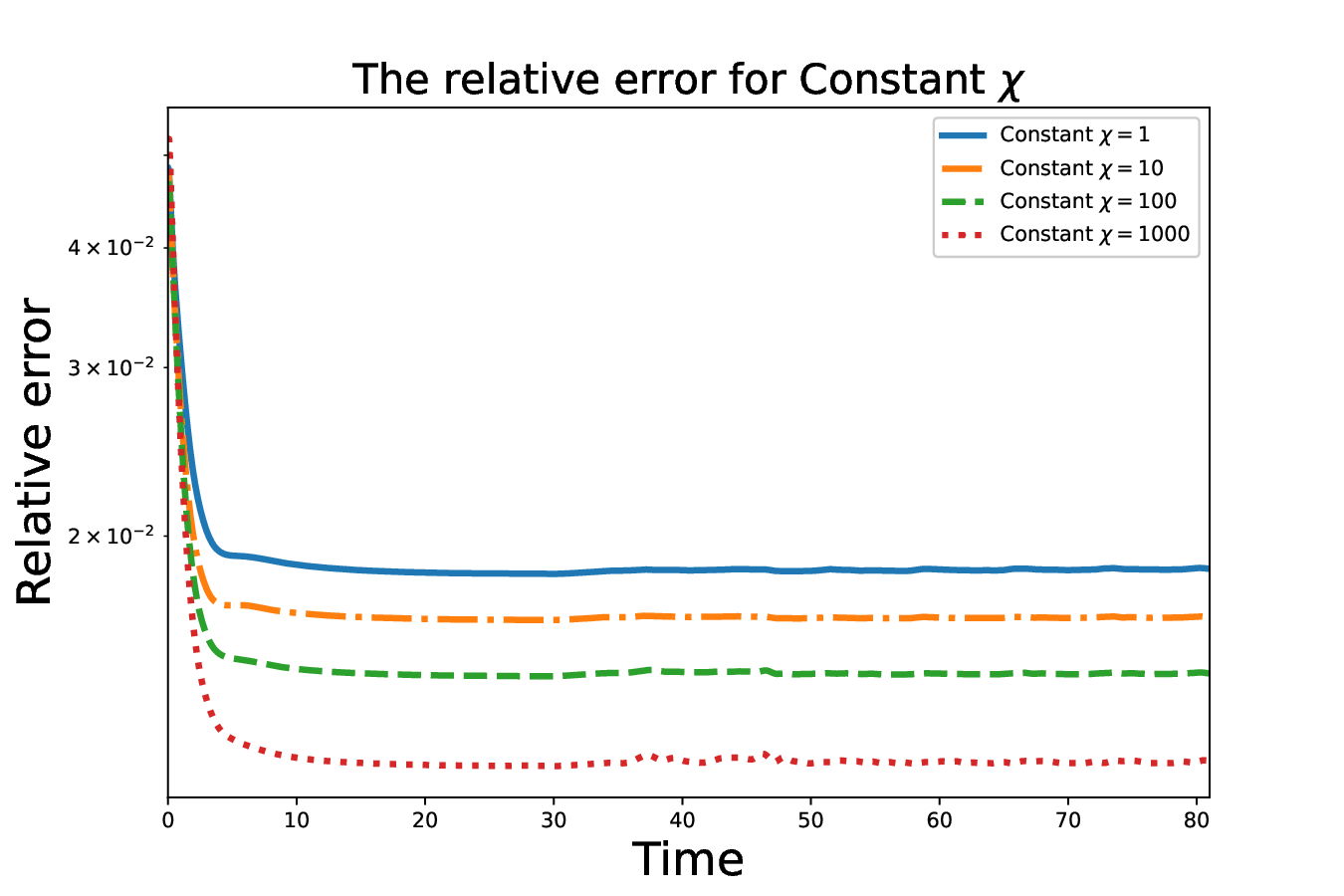}
    \caption{Constant $\chi$: relative errors decrease with larger $\chi$ value.}
    \label{fig: Leo-const-rela-err}
\end{figure}

\begin{figure}[H]
     \centering
    \includegraphics[width=0.8\linewidth]{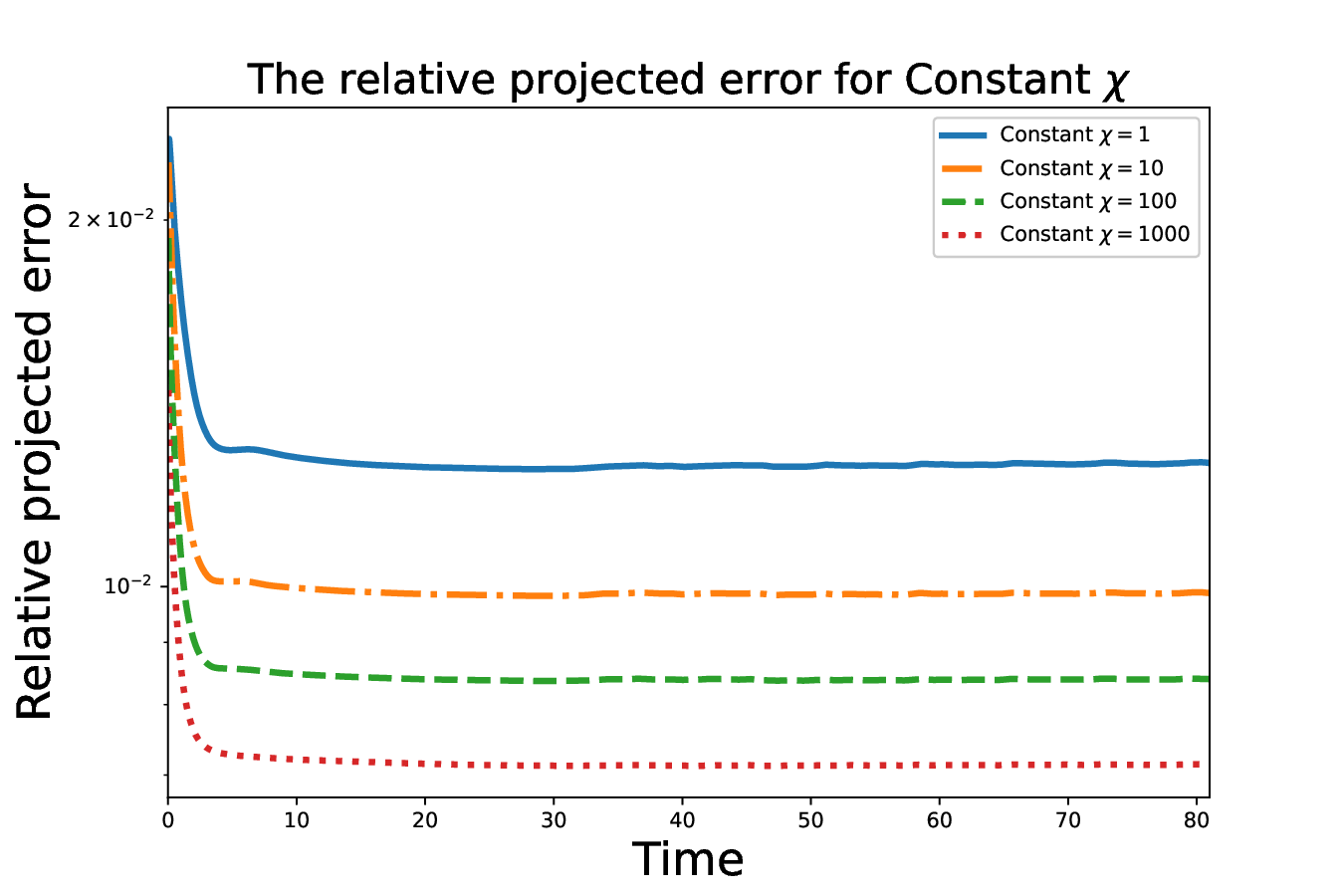}
    \caption{Constant $\chi$: projection errors decrease with larger $\chi$ values.}
    \label{fig: Leo-const-proj-err}
\end{figure}

Since this test is a through-flow problem, it interrogates errors over 1 through-flow time, not over longer times. Still, it is a complex flow with many interesting features and an accepted test problem. 
We compare Algorithm $1$, Algorithm $2$ and the constant $\chi$ with different initial $\chi$ values, where $\chi = 1, 10, 100$, and $1000$. We set $\text{Factor} =1$ and $Tol =0.2$. We define the relative error $\frac{\|u-v\|}{\|u\|}$, and relative projection error $\frac{\|I_H(u-v)\|}{\|u\|
}$.  In Figure \ref{fig: Leo-algo1-rela-err}--\ref{fig: Leo-const-proj-err}, both the relative error and the relative projection error decrease with larger $\chi$ values for Algorithm 1 and the constant $\chi$. Algorithm 1 performs better when the initial $\chi$ value is small because it effectively adapts $\chi$ to reduce the projection error. On the other hand, Algorithm 2 is robust to initial $\chi$ values, adjusting $\chi$ to obtain the optimal value for smaller errors at each time step. We observe that the projection error is smaller than the error. Furthermore, the $\chi$ value for Algorithm 1 is smaller than Algorithm 2 in the first turnover time, and it is bigger than Algorithm 2 in long times. We obtain smaller relative error and projection errors for Algorithm 2 than for Algorithm 1 and the constant $\chi$. We would like to point out here that if $u$ and $v$ are solved on the same mesh (as in, e.g., \cite{leotest}) instead of fine mesh coarse mesh convention as has been done here, the error is negligible, $\mathcal{O}$(machine epsilon).

\begin{figure}[H]
     \centering         \includegraphics[width=0.8\linewidth]{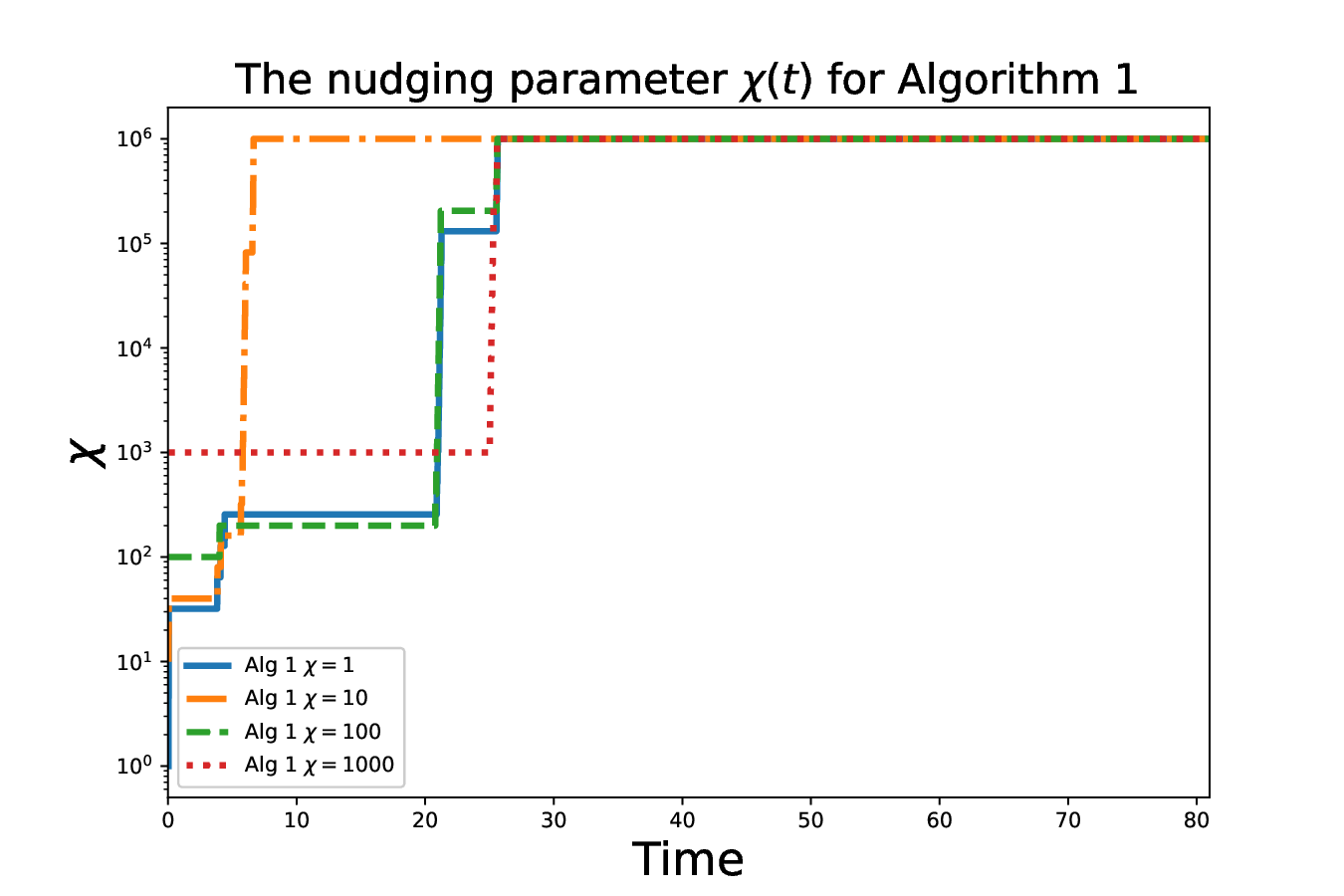}
    \caption{Evolution of $\chi$ values for Algorithm 1. In the first turnover time, from $t=0$ to $t=27$, Algorithm 1 increases the $\chi$ values for smaller initial $\chi$ until it reaches the maximum value of $\chi_{\max} = 10^{6}$.}
    \label{fig: Leo-chi-algo1}
\end{figure}

\section{Conclusion}
Adaptive selection of the nudging parameter $\chi$ improves performance. It does not, however, eliminate the necessity of the $H$-condition related to the density of observations. This means revising the standard formulation of nudging to improve the $H$-condition is an important open problem. Other critical open problems include the effect of time delays in the nudging control term and the effectiveness of nudging for correction of model errors.

\section*{Acknowledgments}
This research of W. Layton, R. Fang, and F. Siddiqua was supported in part by the
NSF under grants DMS 2110379 and 2410893. The author Aytekin \c{C}{\i}b{\i}k was partially supported by TUBITAK with BIDEB-2219 grant.

\bibliography{mybib}
\end{sloppypar}
\end{document}